\documentclass[final,leqno,onefignum,onetabnum]{siamltex1213}
\usepackage{hyperref}
\usepackage{graphicx,epsfig,color}
\usepackage{bm}
\usepackage{booktabs}
\usepackage[notref,notcite]{showkeys}
\usepackage{amsmath,cases}

\usepackage{amssymb}%,pxfonts}

\newtheorem{remark}[theorem]{Remark}
\renewcommand{\theequation}{\arabic{section}.\arabic{equation}}
\numberwithin{equation}{section}

\input{siam.sty}

\title{Preconditioning fractional spectral collocation%\thanks{Version date: 2015.10.13}
}

\author{Kui Du\thanks{%Corresponding author. 
School of Mathematical Sciences and Fujian Provincial Key Laboratory of Mathematical Modeling and High-Performance Scientific Computation, Xiamen University, Xiamen 361005, China ({kuidu@xmu.edu.cn}). The research of this author was supported by the National Natural Science Foundation of China (No.11201392 and No.91430213), the Doctoral Fund of Ministry of Education of China (No.20120121120020), and the Fundamental Research Funds for the Central Universities (No.2013121003).} 
}
\begin{document}

\maketitle

\begin{abstract} Fractional spectral collocation (FSC) method based on fractional Lagrange interpolation has recently been proposed to solve fractional differential equations. Numerical experiments show that the  linear systems in FSC become extremely ill-conditioned as the number of collocation points increases. By introducing suitable fractional Birkhoff interpolation problems, we present fractional integration preconditioning matrices for the ill-conditioned linear systems in FSC. The condition numbers of the resulting linear systems are independent of the number of collocation points. Numerical examples are given.
\end{abstract}

\begin{keywords} 
Fractional Lagrange interpolation, fractional Birkhoff interpolation, fractional spectral collocation, preconditioning
\end{keywords}

\begin{AMS} 65L60, 41A05, 41A10
\end{AMS}

\pagestyle{myheadings}
\thispagestyle{plain}
\markboth{Kui Du}{Preconditioning fractional spectral collocation}

\section{Introduction} Fractional spectral collocation (FSC) methods \cite{zayernouri2014fract,zayernouri2015fract,fatone2015optim} based on fractional Lagrange interpolation have recently been proposed to solve fractional differential equations. By a spectral theory developed in \cite{zayernouri2013fract} for fractional Sturm-Liouville eigenproblems, the corresponding fractional differential matrices can be obtained with ease. However, numerical experiments show that the involved linear systems become extremely ill-conditioned as the number of collocation points increases. Typically, the condition number behaves like $\mcalo(N^{2\nu})$, where $N$ is the number of collocation points and $\nu$ is the order of the leading fractional term. Efficient preconditioners are highly required when solving the linear systems by an iterative method. 

Recently, Wang, Samson, and Zhao \cite{wang2014well} proposed a well-conditioned collocation method to solve linear differential equations with various types of boundary conditions. By introducing a suitable Birkhoff interpolation problem, they constructed a pseudospectral integration preconditioning matrix, which is the exact inverse of the pseudospectral discretization matrix of the $n$th-order derivative operator together with $n$ boundary conditions. Essentially, the linear system in the well-conditioned collocation method \cite{wang2014well} is the one obtained by right preconditioning the original linear system; see \cite{du2015prersc}. By introducing suitable fractional Birkhoff interpolation problems and employing the same techniques in \cite{wang2014well}, Jiao, Wang, and Huang \cite{jiao2015well} proposed fractional integration preconditioning matrices for linear systems in fractional collocation methods base on Lagrange interpolation. In the Riemann-Liouville case, it is necessary to modify the fractional derivative operator in order to absorb  singular fractional factors (see \cite[\S 3]{jiao2015well}).
In this paper, we extend the Birkhoff interpolation preconditioning techniques in \cite{wang2014well,jiao2015well} to the fractional spectral collocation methods \cite{zayernouri2014fract,zayernouri2015fract,fatone2015optim}
based on fractional Lagrange interpolation. Unlike that in \cite{jiao2015well}, there are no singular fractional factors in the Riemann-Liouville case. Numerical experiments show that the condition number of the resulting linear system is independent of the number of collocation points.

The rest of the paper is organized as follows. In \S 2, we review several topics required in the following sections. In \S 3, we introduce fractional Birkhoff interpolation problems and the corresponding fractional integration matrices. In \S 4, we present the preconditioning fractional spectral collocation method. Numerical examples are also reported. We present brief concluding remarks in \S 5.

\section{Preliminaries} 

\subsection{Fractional derivatives}
The definitions of fractional derivatives of order $\nu\in(n-1,n), n\in \mbbn$, on the interval $[-1,1]$ are as follows \cite{kilbas2006theor}:

\begin{itemize}

\item Left-sided Riemann-Liouville fractional derivative: 
$$\rlfdl u(x)=\frac{1}{\Gamma(n-\nu)}\frac{\rmd^n}{\rmd x^n}\int^x_{-1}\frac{u(\xi)}{(x-\xi)^{\nu-n+1}}\rmd\xi,$$

\item Right-sided Riemann-Liouville fractional derivative: 
$$\rlfdr u(x)=\frac{(-1)^n}{\Gamma(n-\nu)}\frac{\rmd^n}{\rmd x^n}\int^1_x\frac{u(\xi)}{(\xi-x)^{\nu-n+1}}\rmd\xi,$$

\item Left-sided Caputo fractional derivative: 
$$\cfdl u(x)=\frac{1}{\Gamma(n-\nu)}\int^x_{-1}\frac{u^{(n)}(\xi)}{(x-\xi)^{\nu-n+1}}\rmd\xi,$$

\item Right-sided Caputo fractional derivative: 
$$\cfdr u(x)=\frac{(-1)^n}{\Gamma(n-\nu)}\int^1_x\frac{u^{(n)}(\xi)}{(\xi-x)^{\nu-n+1}}\rmd\xi.$$

\end{itemize} 

By the definitions of fractional derivatives,  we have 
\beq\label{rlcl} \rlfdl u(x)=\sum_{i=0}^{n-1}\frac{u^{(i)}(-1)}{\Gamma(1+i-\nu)}(x+1)^{i-\nu}+ \cfdl u(x),\eeq  and
\beq\label{rlcr} \rlfdr u(x)=\sum_{i=0}^{n-1}\frac{(-1)^iu^{(i)}(1)}{\Gamma(1+i-\nu)}(1-x)^{i-\nu}+ \cfdr u(x).\eeq 
Therefore, \beqs  \rlfdl u(x)=\cfdl u(x),\quad {\rm if }\quad u^{(i)}(-1)=0,\quad i=0,1,\ldots,n-1,\eeqs and \beqs \rlfdr u(x)=\cfdr u(x),\quad {\rm if }\quad u^{(i)}(1)=0,\quad i=0,1,\ldots,n-1.\eeqs

In this paper, we mainly deal with the left-sided Riemann-Liouville fractional problems with homogeneous boundary/initial conditions. By a simple change of variables, (\ref{rlcl}) and (\ref{rlcr}), the extension to other fractional problems is easy. 

\subsection{Fractional Lagrange interpolation} 

Throughout the paper, let $\{x_j\}_{j=1}^N$ be a set of distinct points satisfying \beq\label{iptsx} -1<x_1<\cdots <x_{N-1}<x_N\leq 1.\eeq 
 %The {\it barycentric resampling matrix}, ${\bf P}^{{\bf x}\mapsto {\bf y}}\in{\mbbr^{M\times N}}$, which interpolates between the points $\{x_j\}_{j=1}^N$ and $\{y_j\}_{j=1}^M$, is defined as \beqs{\bf P}^{{\bf x}\mapsto {\bf y}}=[p^{{\bf x}\mapsto {\bf y}}_{ij}]_{i=1,j=1}^{M,N},\eeqs where \beqs p^{{\bf x}\mapsto {\bf y}}_{ij}=\l\{\begin{array}{ll} \dsp\frac{w_j}{y_i-x_j}\l(\sum_{l=1}^N\dsp\frac{w_l}{y_i-x_l}\r)^{-1}, & y_i\neq x_j, \\ 1, & y_i=x_j, \end{array} \r.\eeqs  and $\{w_j\}_{j=1}^N$ are the barycentric weights associated to the points $\{x_j\}_{j=1}^N$ defined as \beq\label{bweight}w_j=\prod_{n=1,n\neq j}^N(x_j-x_n)^{-1}, \qquad j=1,\ldots,N.\eeq
Given $\mu\in(0,1)$, the {\it $\mu$-fractional Lagrange interpolation basis} associated with the points $\{x_j\}_{j=1}^N$ is defined as \beq\label{mufl} \ell_{j}^\mu(x)=\l(\frac{x+1}{x_j+1}\r)^{\mu}\prod_{n=1,n\neq j}^N\frac{x-x_n}{x_j-x_n},\qquad j=1,\ldots,N.\eeq 
For a function $u(x)$ with $u(-1)=0$, the $\mu$-fractional Lagrange interpolant $u_N(x)$ of $u(x)$ takes the form $$u_N(x)=\sum_{j=1}^Nu(x_j)\ell_{j}^\mu(x).$$ 

%Given $\{x_j\}_{j=1}^N$ as the collocation points, we have the {\it $\nu$-fractional $N\times N$ differentiation matrix}: \beqs\label{psdm1} {\bf D}^{(\nu)}_{{\bf x}\mapsto{\bf x}}=\l[\rlfdl \ell_{j}^{\mu}(x_i)\r]_{i,j=1}^N. \eeqs
% Let $\{y_j\}_{j=1}^M$ be a set of distinct points satisfying \beq\label{iptsy}-1< y_1< \cdots <y_{M-1}<y_M\leq 1.\eeq Define the {\it $\nu$-fractional $M\times N$ differentiation matrix}: \beqs\label{psdm2} {\bf D}^{(\nu)}_{{\bf x}\mapsto {\bf y}}=\l[\rlfdl \ell_{j}^{\mu}(y_i)\r]_{i=1,j=1}^{M,N}.  \eeqs 
%It is easy to show that ${\bf D}^{(\nu)}_{{\bf x}\mapsto {\bf y}}={\bf P}^{{\bf x}\mapsto {\bf y}}{\bf D}^{(\nu)}_{{\bf x}\mapsto{\bf x}}.$

\subsection{Computations of ${^{RL}_{-1}\mcald^{\mu}_x} \ell_{j}^{\mu}(x)$ and ${^{RL}_{-1}\mcald^{1+\mu}_x} \ell_{j}^{\mu}(x)$ with $\mu\in(0,1)$} 

Note that $\ell_{j}^\mu(x)$, $j=1,\ldots,N,$ can be represented exactly as \beq\label{hexpan} \ell_{j}^\mu(x)=\l(\frac{x+1}{x_j+1}\r)^{\mu}\prod_{n=1,n\neq j}^N\frac{x-x_n}{x_j-x_n}=(x+1)^{\mu}\sum_{n=1}^N\alpha_{nj} P^{(-\mu,\mu)}_{n-1}(x),\eeq where $P^{(\alpha,\beta)}_n(x)$ denote the standard Jacobi polynomials. 
The coefficients $\alpha_{nj}$ can be obtained by solving the linear system $$\sum_{n=1}^N (x_i+1)^{\mu}P^{(-\mu,\mu)}_{n-1}(x_i)\alpha_{nj}=\delta_{ij},\quad i=1,\ldots,N.$$ 

\begin{remark}\label{rem1} Let $\{x_j\}_{j=1}^N$ and $\{\omega_j\}_{j=1}^N$ be the Gauss-Jacobi quadrature nodes and weights with the Jacobi polynomial $P_N^{(-\mu,\mu)}(x)$. Then, \beqs \alpha_{nj}=\frac{1}{(x_j+1)^\mu}\frac{(2n-1)(n-1)!(n-1)!}{2\G(n-\mu)\G(n+\mu)}\omega_jP_{n-1}^{(-\mu,\mu)}(x_j).\eeqs \end{remark} 

%\begin{remark} Let $\{x_j^{(-\mu,\mu)}\}_{j=1}^N$ and $\{\omega_j^{(-\mu,\mu)}\}_{j=1}^N$ be the Gauss-Jacobi quadrature nodes and weights with the Jacobi polynomial $P_N^{(-\mu,\mu)}(x)$. Then, \beqs \alpha_{nj}=\frac{1}{(x_j^{(-\mu,\mu)}+1)^\mu}\frac{(2n-1)(n-1)!(n-1)!}{2\G(n-\mu)\G(n+\mu)}\omega_j^{(-\mu,\mu)}P_{n-1}^{(-\mu,\mu)}(x_j^{(-\mu,\mu)}).\eeqs \end{remark} 

%\begin{remark} Let $\{\wt x_j^{(-\mu,\mu)}\}_{j=0}^N$ and $\{\wt \omega_j^{(-\mu,\mu)}\}_{j=0}^N$ be the Gauss-Jacobi-Lobatto quadrature nodes and weights with the Jacobi polynomial $P_N^{(-\mu,\mu)}(x)$. Then, $$\alpha_{nj}=\frac{1}{(\wt x_j^{(-\mu,\mu)}+1)^\mu}\frac{(2n-1)(n-1)!(n-1)!}{2\G(n-\mu)\G(n+\mu)}\l(\wt\omega_j^{(-\mu,\mu)}P_{n-1}^{(-\mu,\mu)}\l(\wt x_j^{(-\mu,\mu)}\r)-\prod_{k=1,k\neq j}^N\frac{-1-\wt x^{(-\mu,\mu)}_k}{\wt x^{(-\mu,\mu)}_j-\wt x^{(-\mu,\mu)}_k}\wt w^{(-\mu,\mu)}_0P^{(-\mu,\mu)}_{n-1}(-1)\r).$$ \end{remark} 

We now compute ${^{RL}_{-1}\mcald^{\mu}_x} \ell_{j}^{\mu}(x)$ and ${^{RL}_{-1}\mcald^{1+\mu}_x} \ell_{j}^{\mu}(x)$. Let $P_n(x)$ denote the Legendre polynomial of order $n$. By (see \cite{zayernouri2013fract}) \beqs {^{RL}_{-1}\mcald^\mu_x}\l((x+1)^\mu P_{n-1}^{(-\mu,\mu)}(x)\r)=\frac{\G(n+\mu)}{\G(n)}P_{n-1}(x)\eeqs and $${^{RL}_{-1}\mcald^{1+\mu}_x} \ell_{j}^{\mu}(x)=\frac{\rmd}{\rmd x}\l({^{RL}_{-1}\mcald^\mu_x}\ell_{j}^{\mu}(x)\r),$$ we have $$ {^{RL}_{-1}\mcald^\mu_x} \ell_{j}^{\mu}(x)=\sum_{n=1}^N\alpha_{nj} \frac{\G(n+\mu)}{\G(n)}P_{n-1}(x)$$ and \beqas  {^{RL}_{-1}\mcald^{1+\mu}_x} \ell_{j}^{\mu}(x)&=&\sum_{n=2}^N\alpha_{nj} \frac{\G(n+\mu)}{\G(n)}P_{n-1}'(x)\\ &=&\sum_{n=2}^N\alpha_{nj} \frac{\G(n+\mu)}{\G(n)}\frac{n}{2}P^{(1,1)}_{n-2}(x).\eeqas

\section{Riemann-Liouville fractional Birkhoff interpolation} 
Let $\mbbp_n$ be the set of all algebraic polynomials of degree at most $n$. Define the space $$\mbbs^\mu_N=(x+1)^\mu\mbbp_{N-1}.$$ %For $\nu=m+\mu\in(m,m+1)$, given points $\{y_j\}_{j=1}^{N-m}$,  the Riemann-Liouville $\nu$-fractional Birkhoff interpolation problem is defined as: find $p(x)\in\mbbs_N^\mu$ such that \beq\label{birkhoff} \l\{\begin{array}{ll}p(-1)= u(-1),& \\ \rlfdl p(y_j)=\rlfdl u(y_j),& j=1,\ldots,N-m,\\ \mcall_i(p)=\mcall_i(u), & i=1,\ldots,m, \end{array}\r.\eeq where each $\mcall_i$ is a linear functional. 
In the following, we consider two special cases.

\subsection{The case $\nu=\mu\in(0,1)$} For a function $u(x)$ with $u(-1)=0$, given $N$ distinct points $\{y_j\}_{j=1}^N$ satisfying $$ -1<y_1<\cdots <y_{N-1}<y_N\leq 1,$$ consider the Riemann-Liouville $\nu$-fractional Birkhoff interpolation problem: 
\beq\label{birkhoff1} {\rm Find}\ \ p(x)\in\mbbs_N^\mu\ \ {\rm such\ that}\ \ \rlfdl p(y_j)=\rlfdl u(y_j),\quad j=1,\ldots,N.\eeq

\begin{theorem}\label{tcase1} The interpolant $u_N^\nu(x)$ for the Riemann-Liouville $\nu$-fractional Birkhoff  problem {\rm(\ref{birkhoff1})} of a function $u(x)$ with $u(-1)=0$ takes the form $$u_N^\nu(x)=\sum_{j=1}^{N}\rlfdl u(y_j)B_j^\nu(x),$$ where \beqs B_j^\nu(x)=(x+1)^\nu\sum_{n=1}^N\wt\alpha_{nj}P_{n-1}^{(-\nu,\nu)}(x),\quad j=1,\ldots,N,\eeqs with $\wt\alpha_{nj}$ satisfying $$\sum_{n=1}^N\wt\alpha_{nj}\frac{\G(n+\nu)}{\G(n)}P_{n-1}(y_i)=\delta_{ij},\quad i=1,\ldots,N.$$ 
\end{theorem}

By $\rlfdl B_j^\nu(y_i)=\delta_{ij}$, the proof of Theorem \ref{tcase1} is straightforward. 

\begin{remark}\label{rem2} Let $\{y_j\}_{j=1}^{N}$ and $\{\omega_j\}_{j=1}^{N}$ be the Gauss-Legendre quadrature nodes and weights with the Legendre polynomial $P_N(x)$. Then, %\beqas && B_0^\nu(x)=1, \\ && B_j^\nu(x)=\sum_{n=1}^N\wt\alpha_{nj}(x+1)^\nu P_{n-1}^{(-\nu,\nu)}(x),\quad j=1,\ldots,N,\eeqas and $\wt\alpha_{nj}$ satisfy $$\sum_{n=1}^N\wt\alpha_{nj}\frac{\G(n+\nu)}{\G(n)}P_{n-1}(x)=\ell_{j}(x)=\prod_{n=1,n\neq j}^N\frac{x-y_n}{y_j-y_n}.$$ 
$$\wt\alpha_{nj}=\frac{2n-1}{2}\frac{\G(n)}{\G(n+\nu)}\omega_jP_{n-1}(y_j).$$\end{remark}

Define the matrices $${\bf D}^{(\nu)}_{{\bf x}\mapsto {\bf y}}=\l[\rlfdl \ell_{j}^{\mu}(y_i)\r]_{i,j=1}^{N},\qquad {\bf B}_{{\bf y}\mapsto {\bf x}}^{(-\nu)}=\l[B_j^\nu(x_i)\r]_{i,j=1}^N.$$ It is easy to show that \beq\label{case1}{\bf D}^{(\nu)}_{{\bf x}\mapsto {\bf y}}{\bf B}_{{\bf y}\mapsto {\bf x}}^{(-\nu)}={\bf I}_N.\eeq

\subsection{The case $\nu=1+\mu\in(1,2)$} For a function $u(x)$ with $u(\pm1)=0$, given $N-1$ distinct points $\{y_j\}_{j=1}^{N-1}$ satisfying $$ -1<y_1<\cdots <y_{N-1}<1,$$ consider the Riemann-Liouville $\nu$-fractional Birkhoff interpolation problem: 
\beq\label{birkhoff2} {\rm Find}\ \ p(x)\in\mbbs_N^\mu\ \ {\rm such\ that}\ \ \dsp\l\{\begin{array}{ll}\dsp p(1)=0, & \\ \dsp\rlfdl p(y_j)=\rlfdl u(y_j),& j=1,\ldots,N-1.\end{array}\r.\eeq

\begin{theorem}\label{tcase2} The interpolant $u_N^\nu(x)$ for the Riemann-Liouville $\nu$-fractional Birkhoff  problem {\rm(\ref{birkhoff2})} of a function $u(x)$ with $u(\pm 1)=0$ takes the form $$u_N^\nu(x)=\sum_{j=1}^{N-1}\rlfdl u(y_j)B_j^\nu(x),$$ where  
$$B_j^\nu(x)=(x+1)^\mu\sum_{n=1}^{N-1}\wt\beta_{nj}\l(P_{n}^{(-\mu,\mu)}(x)-P_{n}^{(-\mu,\mu)}(1)\r),\quad j=1,\ldots,N-1,$$
with $\mu=\nu-1$ and $\wt\beta_{nj}$ satisfying \beqs \sum_{n=1}^{N-1}\wt\beta_{nj}\frac{\G(n+1+\mu)}{\G(n+1)}\frac{n+1}{2}P_{n-1}^{(1,1)}(y_i)=\delta_{ij},\quad i=1,\ldots,N-1. \eeqs
\end{theorem}

By $\rlfdl B_j^\nu(y_i)=\delta_{ij}$, the proof of Theorem \ref{tcase2} is straightforward. 

\begin{remark}\label{rem3} Let $\{y_j\}_{j=1}^{N-1}$ and $\{\omega_j\}_{j=1}^{N}$ be the Gauss-Jacobi quadrature nodes and weights with the Jacobi polynomial $P_{N-1}^{(1,1)}(x)$. Then, 
\beqs\wt\beta_{nj}=\frac{2n+1}{4n}\frac{\G(n+1)}{\G(n+1+\mu)}\omega_jP_{n-1}^{(1,1)}(y_j).\eeqs\end{remark}

%\begin{remark} Let $\{y_j^{(1,1)}\}_{j=1}^{N-1}$ and $\{\omega_j^{(1,1)}\}_{j=1}^{N}$ be the Gauss-Jacobi quadrature nodes and weights with the Jacobi polynomial $P_{N-1}^{(1,1)}(x)$. Then, \beqs\wt\beta_{nj}=\frac{2n+1}{4n}\frac{\G(n+1)}{\G(n+1+\mu)}\omega_j^{(1,1)}P_{n-1}^{(1,1)}(y_j^{(1,1)}).\eeqs\end{remark}

In this subsection, let $x_N=1$. Define the matrices $$\wt{\bf D}^{(\nu)}_{{\bf x}\mapsto {\bf y}}=\l[\rlfdl \ell_{j}^{\nu}(y_i)\r]_{i,j=1}^{N-1}, \qquad \wt{\bf B}_{{\bf y}\mapsto {\bf x}}^{(-\nu)}=\l[B_j^\nu(x_i)\r]_{i,j=1}^{N-1}.$$ It is easy to show that \beq\label{case2}\wt{\bf D}^{(\nu)}_{{\bf x}\mapsto {\bf y}}\wt{\bf B}_{{\bf y}\mapsto {\bf x}}^{(-\nu)}={\bf I}_{N-1}.\eeq

\section{Preconditioning fractional spectral collocation (PFSC)} 
In this section, we use two examples to introduce the preconditioning scheme.

\subsection{An initial value problem}
Consider the fractional differential equation of the form \beq\label{fde1}\rlfdl u(x)+a(x)u(x)=f(x), \qquad \nu\in(0,1);\qquad u(-1)=0.\eeq The fractional spectral collocation scheme leads to the following linear system 
\beq\label{fsc1} \l({\bf D}_{{\bf x}\mapsto{\bf y}}^{(\nu)}+\diag\{{\bf a}\}{\bf D}_{{\bf x}\mapsto{\bf y}}^{(0)}\r){\bf u}={\bf f}, \eeq where $${\bf D}_{{\bf x}\mapsto{\bf y}}^{(0)}=\l[\ell_{j}^\nu(y_i)\r]_{i,j=1}^N,$$ and \beqas &&{\bf a}=\l[\begin{array}{cccc}a(y_1)& a(y_2)& \cdots &a(y_{N})\end{array}\r ]^\rmt,\\ &&{\bf f}=\l[\begin{array}{cccc}f(y_1)& f(y_2)& \cdots &f(y_{N})\end{array}\r ]^\rmt.\eeqas The unknown vector ${\bf u}$ is an approximation of the vector of the exact solution $u(x)$ at the points $\{x_j\}_{j=1}^{N}$, i.e., $$\l[\begin{array}{cccc}u(x_1)& u(x_2)& \cdots &u(x_{N})\end{array}\r ]^\rmt.$$

Consider the matrix ${\bf B}_{{\bf y}\mapsto {\bf x}}^{(-\nu)}$ as a right preconditioner for the linear system (\ref{fsc1}). By (\ref{case1}), we have the right preconditioned linear system \beq\label{pfsc1}\l({\bf I}_N+\diag\{{\bf a}\}{\bf D}_{{\bf x}\mapsto{\bf y}}^{(0)}{\bf B}_{{\bf y}\mapsto {\bf x}}^{(-\nu)}\r){\bf v}={\bf f}.\eeq   
It is easy to show that $${\bf D}_{{\bf x}\mapsto {\bf y}}^{(0)}{\bf B}_{{\bf y}\mapsto {\bf x}}^{(-\nu)}={\bf B}_{{\bf y}\mapsto {\bf y}}^{(0-\nu)},$$ where $${\bf B}_{{\bf y}\mapsto {\bf y}}^{(0-\nu)}=\l[B_j^{\nu}(y_i)\r]_{i,j=1}^N.$$ Then, the equation (\ref{pfsc1}) reduces to \beq\label{pfsc1t}\l({\bf I}_N+\diag\{{\bf a}\}{\bf B}_{{\bf y}\mapsto {\bf y}}^{(0-\nu)}\r){\bf v}={\bf f}.\eeq
After solving (\ref{pfsc1t}), we obtain ${\bf u}$ by ${\bf u}={\bf B}_{{\bf y}\mapsto {\bf x}}^{(-\nu)}{\bf v}.$

\subsubsection*{Example 1} We consider the fractional differential equation (\ref{fde1}) with $$\nu=0.8,\qquad a(x)=2+\sin(25x).$$ The function $f(x)$ is chosen such that the exact solution of (\ref{fde1}) is $$u(x)=\rme^{x+1}-1+(x+1)^{46/7}.$$ 

\begin{figure}[!htpb]
\centerline{\epsfig{figure=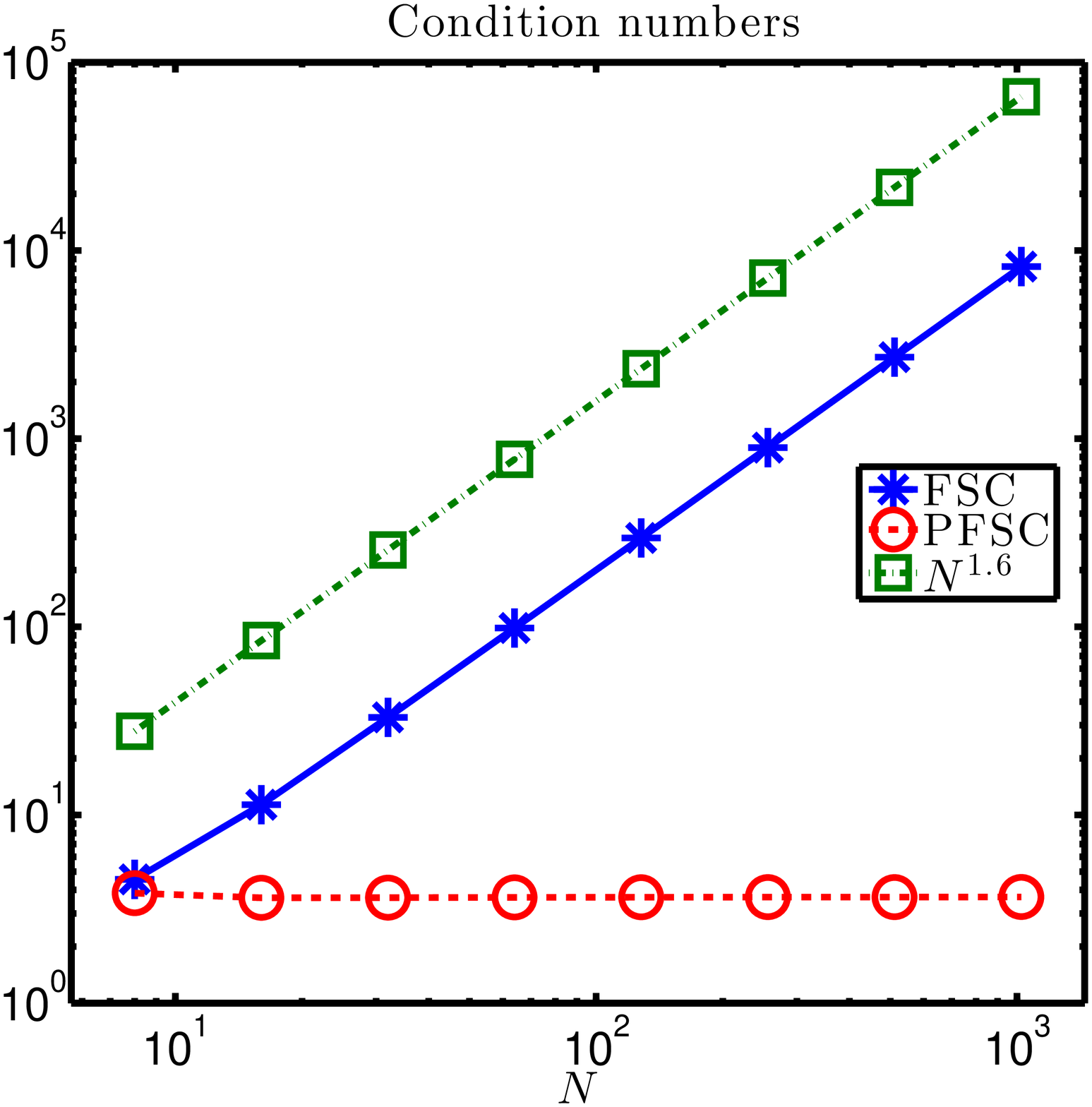,height=1.65in}\hspace{-4mm}\epsfig{figure=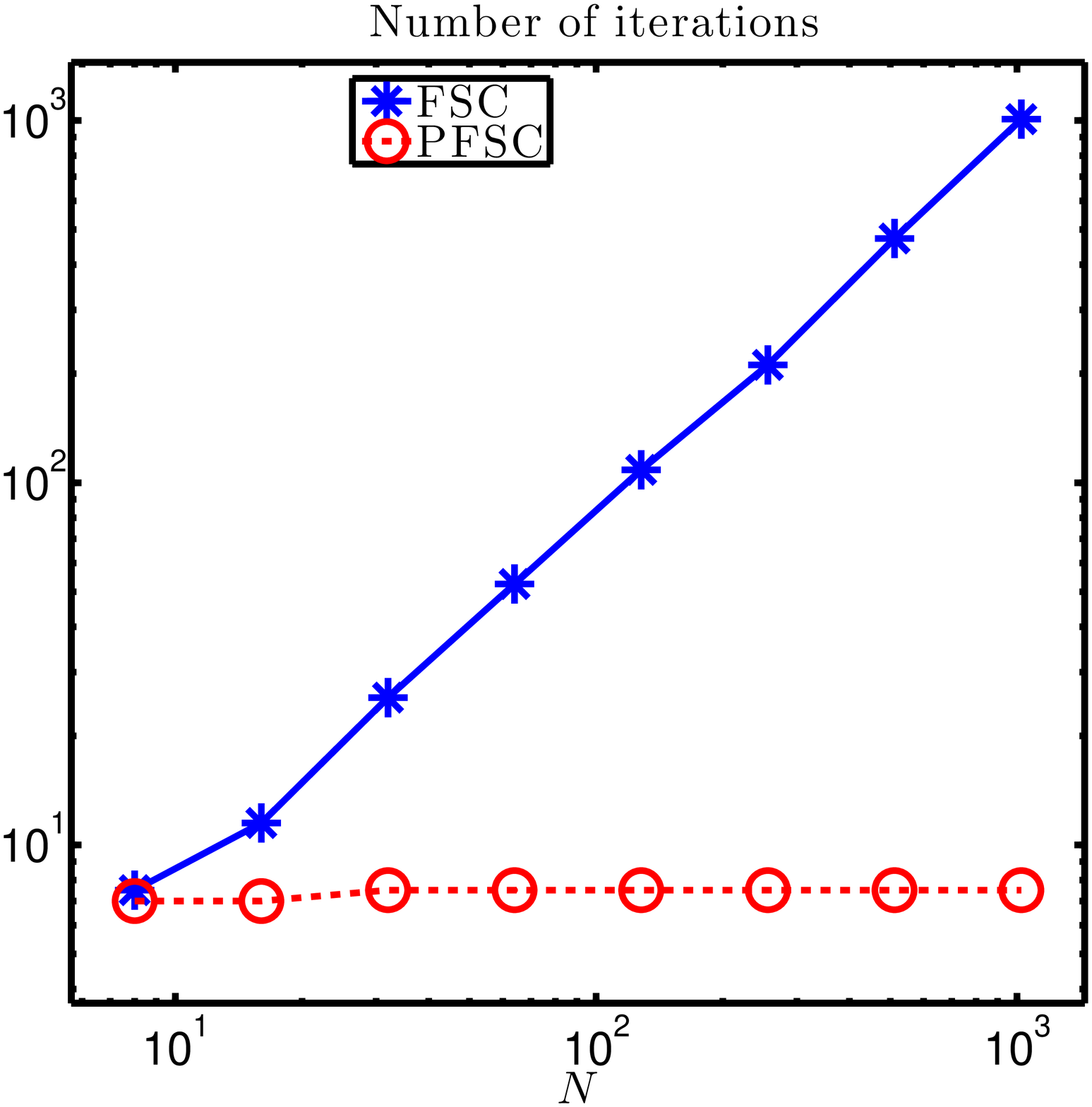,height=1.65in}\hspace{-4mm}\epsfig{figure=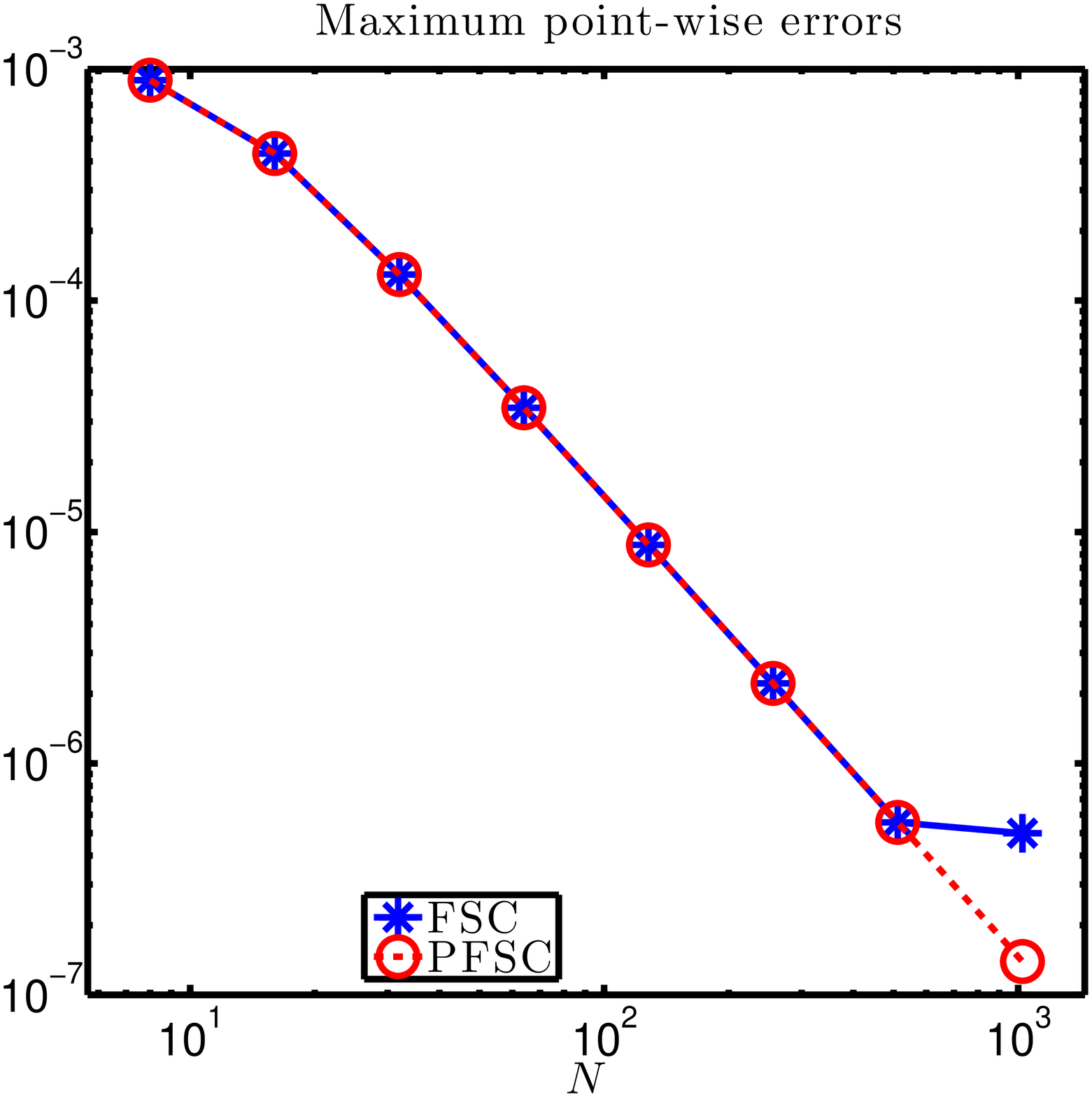,height=1.65in}}
%\centerline{\footnotesize\hspace{3mm}(a)\hspace{3.3cm} (b)\hspace{3.3cm} (c)}
\caption{Comparison of condition numbers {\rm (left)}, number of iterations {\rm (middle)}, and maximum point-wise errors {\rm (right)}  for Example {\rm 1}.}\label{e1f1}
\end{figure}

Let $\{x_j\}_{j=1}^N$ be the Gauss-Jacobi  points as in Remark \ref{rem1} and $\{y_j\}_{j=1}^N$ be the Gauss-Legendre points as in Remark \ref{rem2}. We compare condition numbers, number of iterations (using BiCGSTAB in Matlab with TOL$=10^{-9}$) and maximum point-wise errors of FSC and PFSC (see Figure \ref{e1f1}). Observe from Figure \ref{e1f1} (left) that the condition number of FSC behaves like $\mcalo(N^{1.6})$, while that of PFSC scheme remains a constant even for $N$ up to $1024$. As a result, PFSC scheme only requires about 7 iterations to converge (see  Figure \ref{e1f1} (middle)), while the usual FSC scheme requires much more iterations with a degradation of accuracy as depicted in Figure \ref{e1f1} (right). 

%\begin{table}[htp]
%\caption{Norms of the matrix ${\bf B}_{{\bf y}\mapsto {\bf y}}^{(0-\nu)}$ in Example {\rm 1}.}
%\label{e1t1}
%\begin{center} \footnotesize
%\begin{tabular}{c|c|c|c} \toprule
% $N$& $\|{\bf B}_{{\bf y}\mapsto {\bf y}}^{(0-\nu)}\|_1$ & $\|{\bf B}_{{\bf y}\mapsto {\bf y}}^{(0-\nu)}\|_2$ & $\|{\bf B}_{{\bf y}\mapsto {\bf y}}^{(0-\nu)}\|_\infty$\\ \hline
%  8& 1.7379& 1.3472& 1.8396\\
% 16& 1.7716& 1.3671& 1.8614\\
% 32& 1.7735& 1.3774& 1.8673\\
% 64& 1.7657& 1.3826& 1.8689\\
%128& 1.7569& 1.3852& 1.8692\\
%256& 1.7497& 1.3865& 1.8693\\ \bottomrule
%\end{tabular}
%\end{center}
%\end{table} 

\subsection{A boundary value problem}
Consider the fractional differential equation of the form \beq\label{fde2}\rlfdl u(x)+a(x)u'(x)+b(x)u(x)=f(x), \quad \nu=1+\mu\in(1,2);\quad u(\pm 1)=0.\eeq The fractional spectral collocation method leads to the following linear system 
\beq\label{fsc2} \l(\wt{\bf D}_{{\bf x}\mapsto{\bf y}}^{(\nu)}+\diag\{{\bf a}\}\wt{\bf D}_{{\bf x}\mapsto{\bf y}}^{(1)}+\diag\{{\bf b}\}\wt{\bf D}_{{\bf x}\mapsto{\bf y}}^{(0)}\r){\bf u}={\bf f}, \eeq where $$\wt{\bf D}_{{\bf x}\mapsto{\bf y}}^{(1)}=\l[\frac{\rmd}{\rmd x}\l(\ell^\mu_j(x)\r)\Big|_{x=y_i}\r]_{i,j=1}^{N-1},\qquad \wt{\bf D}_{{\bf x}\mapsto{\bf y}}^{(0)}=\l[\ell^\mu_j(y_i)\r]_{i,j=1}^{N-1},$$ and \beqas &&{\bf a}=\l[\begin{array}{cccc}a(y_1)& a(y_2)& \cdots &a(y_{N-1})\end{array}\r ]^\rmt,\\ &&{\bf b}=\l[\begin{array}{cccc}b(y_1)& b(y_2)& \cdots &b(y_{N-1})\end{array}\r ]^\rmt,\\  &&{\bf f}=\l[\begin{array}{cccc}f(y_1)& f(y_2)& \cdots &f(y_{N-1})\end{array}\r ]^\rmt.\eeqas The unknown vector ${\bf u}$ is an approximation of the vector of the exact solution $u(x)$ at the points $\{x_j\}_{j=1}^{N-1}$, i.e., $$\l[\begin{array}{cccc}u(x_1)& u(x_2)& \cdots &u(x_{N-1})\end{array}\r ]^\rmt.$$

Consider the matrix $\wt{\bf B}_{{\bf y}\mapsto {\bf x}}^{(-\nu)}$ as a right preconditioner for the linear system (\ref{fsc2}). By (\ref{case2}), we have the right preconditioned linear system \beq\label{pfsc2}\l({\bf I}_{N-1}+\diag\{{\bf a}\}\wt{\bf D}_{{\bf x}\mapsto{\bf y}}^{(1)}\wt{\bf B}_{{\bf y}\mapsto {\bf x}}^{(-\nu)}+\diag\{{\bf b}\}\wt{\bf D}_{{\bf x}\mapsto{\bf y}}^{(0)}\wt{\bf B}_{{\bf y}\mapsto {\bf x}}^{(-\nu)}\r){\bf v}={\bf f}.\eeq It is easy to show that $$\wt {\bf D}_{{\bf x}\mapsto {\bf y}}^{(1)}\wt{\bf B}_{{\bf y}\mapsto {\bf x}}^{(-\nu)}=\wt{\bf B}_{{\bf y}\mapsto {\bf y}}^{(1-\nu)},\qquad \wt {\bf D}_{{\bf x}\mapsto {\bf y}}^{(0)}\wt{\bf B}_{{\bf y}\mapsto {\bf x}}^{(-\nu)}=\wt{\bf B}_{{\bf y}\mapsto {\bf y}}^{(0-\nu)},$$ where $$\wt{\bf B}_{{\bf y}\mapsto {\bf y}}^{(1-\nu)}=\l[\frac{\rmd}{\rmd x}\l(B_j^{\nu}(x)\r)\Big|_{x=y_i}\r]_{i,j=1}^{N-1},\qquad \wt{\bf B}_{{\bf y}\mapsto {\bf y}}^{(0-\nu)}=\l[B_j^{\nu}(y_i)\r]_{i,j=1}^{N-1}.$$ Then, the equation (\ref{pfsc2}) reduces to \beq\label{pfsc2t}\l({\bf I}_{N-1}+\diag\{{\bf a}\}\wt{\bf B}_{{\bf y}\mapsto {\bf y}}^{(1-\nu)}+\diag\{{\bf b}\}\wt{\bf B}_{{\bf y}\mapsto {\bf y}}^{(0-\nu)}\r){\bf v}={\bf f}.\eeq  After solving (\ref{pfsc2t}), we obtain ${\bf u}$ by ${\bf u}=\wt{\bf B}_{{\bf y}\mapsto {\bf x}}^{(-\nu)}{\bf v}.$

\subsubsection*{Example 2} We consider the fractional differential equation (\ref{fde2}) with $$\nu=1.9,\qquad a(x)=2+\sin(4\pi x),\qquad b(x)=2+\cos x.$$ The function $f(x)$ is chosen such that the exact solution of (\ref{fde2}) is $$u(x)=\rme^{x+1}-x-2-\frac{\rme^2-3}{4}(x+1)^2+(x+1)^{46/7}-2(x+1)^{39/7}.$$ 

\begin{figure}[!htpb]
\centerline{\epsfig{figure=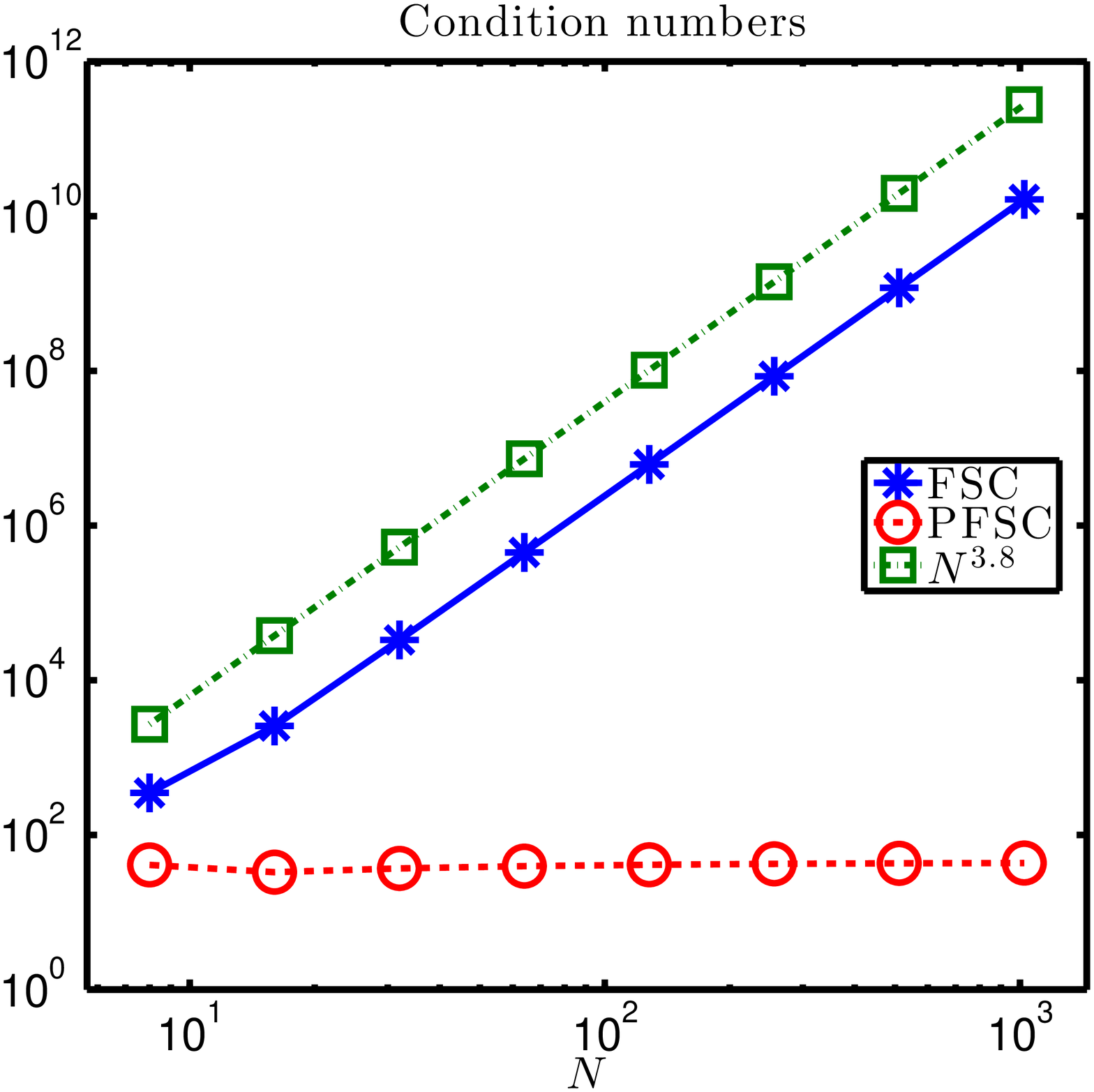,height=1.65in}\hspace{-4mm}\epsfig{figure=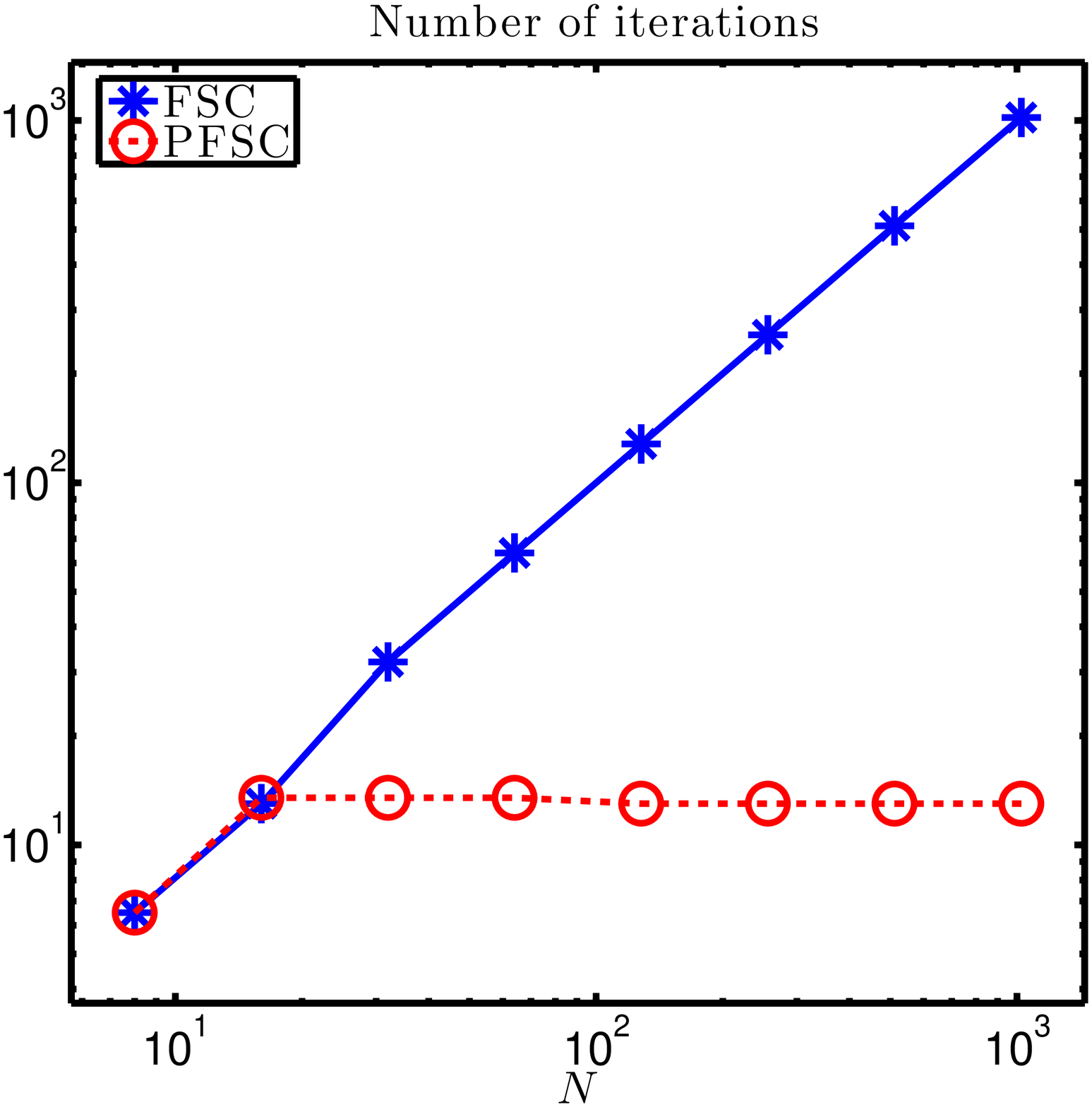,height=1.65in}\hspace{-4mm}\epsfig{figure=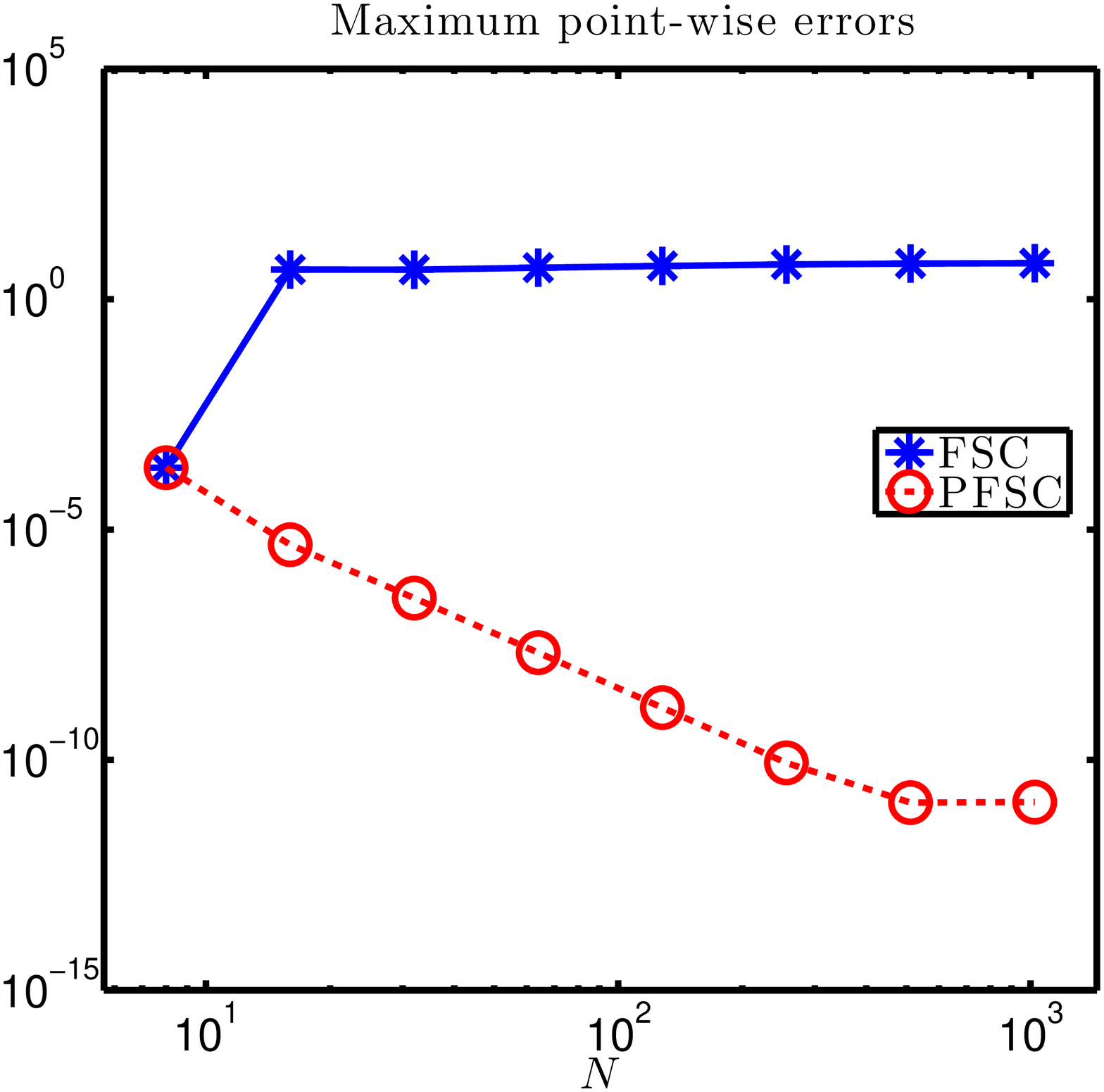,height=1.65in}}
\caption{Comparison of condition numbers {\rm (left)}, number of iterations {\rm (middle)}, and maximum point-wise errors {\rm (right)} for Example {\rm 2}.}\label{e2f1}
\end{figure}

Let $\{x_j\}_{j=0}^N$ be the Chebyshev points of the second kind (also known as Gauss-Chebyshev-Lobatto points) defined as $$x_j=-\cos\frac{j\pi}{N},\qquad  j=0,1,\ldots, N,$$ and $\{y_j\}_{j=1}^{N-1}$ be the Gauss-Jacobi points as in Remark \ref{rem3}. We compare condition numbers, number of iterations (using BiCGSTAB in Matlab with TOL$=10^{-11}$) and maximum point-wise errors of FSC and PFSC (see Figure \ref{e2f1}). Observe from Figure \ref{e2f1} (left) that the condition number of FSC behaves like $\mcalo(N^{3.8})$, while that of PFSC scheme remains a constant even for $N$ up to $1024$. As a result, PFSC scheme only requires about 13 iterations to converge (see  Figure \ref{e2f1} (middle)), while the FSC scheme fails to converge (when $N\geq 16$)  within $N$ iterations as depicted in Figure \ref{e2f1} (right).

%\begin{table}[htp]
%\caption{Norms of the matrices $\wt{\bf B}_{{\bf y}\mapsto {\bf y}}^{(0-\nu)}$ and $\wt{\bf B}_{{\bf y}\mapsto {\bf y}}^{(1-\nu)}$ in Example {\rm 2}.}
%\label{e2t1}
%\begin{center} \footnotesize
%\begin{tabular}{c|c|c|c|c|c|c} \toprule
% $N$& $\|\wt{\bf B}_{{\bf y}\mapsto {\bf y}}^{(0-\nu)}\|_1$ & $\|\wt{\bf B}_{{\bf y}\mapsto {\bf y}}^{(0-\nu)}\|_2$ & $\|\wt{\bf B}_{{\bf y}\mapsto {\bf y}}^{(0-\nu)}\|_\infty$ & $\|\wt{\bf B}_{{\bf y}\mapsto {\bf y}}^{(1-\nu)}\|_1$ & $\|\wt{\bf B}_{{\bf y}\mapsto {\bf y}}^{(1-\nu)}\|_2$ & $\|\wt{\bf B}_{{\bf y}\mapsto {\bf y}}^{(1-\nu)}\|_\infty$\\ \hline
%  8& 0.6227& 0.4493& 0.5472& 1.3373& 0.8266& 1.1086\\
% 16& 0.6260& 0.4494& 0.5472& 1.5414& 0.9112& 1.3741\\
% 32& 0.6268& 0.4494& 0.5477& 1.6580& 0.9638& 1.6073\\
% 64& 0.6276& 0.4494& 0.5487& 1.7227& 0.9965& 1.8535\\
%128& 0.6278& 0.4494& 0.5487& 1.7582& 1.0171& 2.1304\\
%256& 0.6278& 0.4494& 0.5487& 1.7776& 1.0302& 2.4470\\ \bottomrule
%\end{tabular}
%\end{center}
%\end{table} 

\section{Concluding remarks} We numerically show that the Birkhoff interpolation preconditioning techniques in \cite{wang2014well,jiao2015well} are still effective for fractional spectral collocation schemes \cite{zayernouri2014fract,zayernouri2015fract,fatone2015optim} based on fractional Lagrange interpolation. The preconditioned coefficient matrix is a perturbation of the identity matrix. The condition number is independent of the number of collocation points. The preconditioned linear system can be solved by an iterative solver within a few iterations. The application of the preconditioning FSC scheme to multi-term fractional differential equations is straightforward.
%\cite{zayernouri2015temper}

%\bibliographystyle{siam}
%\bibliography{/users/dukui/mywork/common/macbib}

\end{document}